\documentclass[1p]{elsarticle}
\usepackage{amssymb}
\usepackage{amsfonts}

\newtheorem{theorem}{Theorem}
\newtheorem{conjecture}[theorem]{Conjecture}
\newtheorem{corollary}[theorem]{Corollary}

\newtheorem{observation}[theorem]{Observation}

\newtheorem{claim}{Claim}

\newproof{pf}{Proof}

\begin{document}
\title{The 1--2--3 Conjecture almost holds for regular graphs}

\author[agh]{Jakub Przyby{\l}o\fnref{MNiSW}}
\ead{jakubprz@agh.edu.pl} 

\fntext[MNiSW]{This work was partially supported by the Faculty of Applied Mathematics AGH UST statutory tasks within subsidy of Ministry of Science and Higher Education.}

\address[agh]{AGH University of Science and Technology, al. A. Mickiewicza 30, 30-059 Krakow, Poland}

\begin{abstract}
The well-known 1--2--3 Conjecture asserts that the edges of every graph without isolated edges can be weighted with $1$, $2$ and $3$ so that adjacent vertices receive distinct weighted degrees. This is open in general, while it is known to be possible from the weight set $\{1,2,3,4,5\}$. We show that for regular graphs it is sufficient to use weights $1$, $2$, $3$, $4$. Moreover, we prove the conjecture to hold for every $d$-regular graph with $d\geq 10^8$.
\end{abstract}


\maketitle

\section{Introduction}

One of the most basic observations in graph theory implies that there are no antonyms of regular graphs, understood as graphs whose all vertices have pairwise distinct degrees, except the trivial one vertex case.
 Potential alternative definitions of an irregular graph were studied in the paper of Chartrand, Erd\H{o}s and Oellermann~\cite{ChartrandErdosOellermann}.
Chartrand et al.~\cite{Chartrand} on the other hand turned towards measuring the level of irregularity of graphs, rather than defining their irregular representatives. Suppose we admit multiplying edges of a given simple graph $G$, then what is the minimum $k$ such that we may obtain an irregular multigraph (a multigraph with pairwise distinct all degrees) of $G$ via replacing its every edge $e$ by at most $k$ parallel copies of $e$? Such value was called the \emph{irregularity strength} of $G$, see details and exemplary results concerning this graph invariant in~\cite{Aigner,Lazebnik,Faudree,Frieze,KalKarPf,Lehel,MajerskiPrzybylo2,Nierhoff}.
It was investigated in numerous further papers  and gave rise to a wide list of related problems. 
Perhaps the most closely associated with the irregularity strength itself is its local variant, oriented towards 
differentiating degrees of exclusively adjacent vertices. Note that rather than multiplying edges of a given 
graph $G=(V,E)$, we may consider its  \emph{edge $k$-weighting}, i.e. an assignment $\omega:E\to\{1,2,\ldots,k\}$, and instead of focusing on a degree of a vertex $v$ in the corresponding multigraph, consider its so-called \emph{weighted degree} in $G$, defined as $\sigma_\omega(v):=\sum_{u\in N(v)}\omega(uv)$.
If this causes no ambiguities, we also write $\sigma(v)$ instead of $\sigma_\omega(v)$ and call it simply the 
 \emph{sum at $v$}. 
 We say $\omega$ is \emph{vertex-colouring} if $\sigma(u)\neq\sigma(v)$ for every edge $uv\in E$ --
we shall write that $u$ and $v$ are \emph{sum-distinguished} then or that there is \emph{no sum conflict} between $u$ and $v$. 
This concept gained equally considerable attention in the combinatorial community as its precursor largely due to the following intriguing conjecture.
\begin{conjecture}[1--2--3 Conjecture]\label{Conjecture123Conjecture}
Every graph without isolated edges admits a ver\-tex-co\-louring edge $3$-weighting.
\end{conjecture}
This remarkable presumption originates in the paper~\cite{123KLT} of Karo\'nski, {\L}uczak and Thomason, who confirmed it in particular for $3$-colourable graphs.
First general constant upper bound was however showed by Addario-Berry, Dalal, McDiarmid, Reed and Thomason~\cite{Louigi30}, who designed strong and vastly applicative theorems on so-called degree constrained subgraphs (cf. e.g.~\cite{Louigi2}) to prove that every graph without isolated edges admits a vertex-colouring edge $30$-weighting. The same technique was further developed by Addario-Berry,  Dalal and Reed~\cite{Louigi} to decrease $30$ to $16$, and by Wang and Yu~\cite{123with13}, who pushed it further down to $13$. A big break-through was later achieved due to research devoted to a total variant of the same concept, introduced in~\cite{12Conjecture}, and especially thanks to the result of Kalkowski~\cite{Kalkowski12}, generalized later through algebraic approach towards list setting by Wong and Zhu~\cite{WongZhu23Choos}. See also e.g.~\cite{BarGrNiw,PrzybyloWozniakChoos,WongZhuChoos} for other results, concerning in particular list versions of the both problems. 
A modification and development of a surprisingly simple algorithm designed by Kalkowski in~\cite{Kalkowski12} allowed 
Kalkowski, Karo\'nski and Pfender~\cite{KalKarPf_123} to
achieve the best general bound thus far in view of Conjecture~\ref{Conjecture123Conjecture}, 
implying that weights $1,2,3,4,5$ are always sufficient.
It is moreover known that the 1--2--3 Conjecture holds for very dense and large enough graphs, i.e. that there exists a constant $n'$ such that every graph with $n\geq n'$ vertices and minimum degree $\delta(G)>0.99985n$ admits a vertex-colouring edge $3$-weighting, as proved recently by Zhong in~\cite{123dense-Zhong}, and that even just weights 1,2 are asymptotically almost surely sufficient for a random graph (chosen from $G_{n,p}$ for a constant $p\in(0,1)$), see~\cite{Louigi}. On the other hand, it was proved by Dudek and Wajc~\cite{DudekWajc123complexity} that determining whether a particular graph admits a vertex-colouring edge $2$-weighting is NP-complete, while Thomassen, Wu and Zhang~\cite{ThoWuZha} showed that the same problem is polynomial in the family of bipartite graphs.
In this paper we provide two results drawing us very close to a complete solution of the 1--2--3 Conjecture in the case of regular graphs, which apparently might seem most obstinate in its context in view of exactly equal degrees of all their vertices (though obviously regularities within them might be and are an asset while analysing them).

\section{Main Results}

\begin{theorem}\label{1234regTh}
Every $d$-regular graph with $d\geq 2$ admits a vertex-colouring edge $4$-weighting.
\end{theorem}
This was earlier known for $d\leq 3$ by~\cite{123KLT} and possibly for $d=4$ -- see e.g.~\cite{Seamon123survey}.
Recently this was also confirmed for $d=5$ by Bensmail~\cite{Julien5regular123}.
Though our proof of  Theorem~\ref{1234regTh} was obtained independently, its generic idea partly resembles
the one from~\cite{Julien5regular123}, but extends to all regular graphs. Apart from this,  we moreover prove that the 1--2--3 Conjecture holds for regular graphs with sufficiently large degree by showing the following. 
\begin{theorem}\label{123largeregTh}
Every $d$-regular graph with $d\geq 10^8$ admits a vertex-colouring edge $3$-weighting.
\end{theorem}
The proof of this result is completely different and based on the probabilistic method. 
There can be found two common factors of the both approaches though. 
Firstly, both exploit at some point modifications of Kalkowski's algorithm from~\cite{Kalkowski12} in order to get read of a part of possible sum conflicts, but in a different manner -- this is used as one of the main tools in the first proof, and only as a kind of a final cleaning device in the second one. Secondly, in the both approaches we single out a special, usually small subset of vertices, and use the edges between this set an the rest of the vertices to adjust the sums in the graph (in the first proof such a set $I$ is stable, while in the second one -- such set $V_0$ is chosen randomly, and thus usually not stable), all details follow.

We shall apply the following rather standard notation for any given graph $G=(V,E)$, $v\in V$, $E'\subseteq E$ and $V',V''\subseteq V$ where $V'\cap V''=\emptyset$. By $G[V']$ we understand the graph induced by $V'$ in $G$, by $N_{V'}(v)$ -- the set of edges  $uv\in E$ with $u\in V'$, by $d_{V'}(v)$ -- the number of edges $uv\in E$ with $u\in V'$ (i.e., $d_{V'}(v):=|N_{V'}(v)|$), by $d_{E'}(v)$ -- the number of edges in $E'$ incident with $v$, by $E(V')$ -- the set of edges from $E$ with both ends in $V'$, by $E(V',V'')$ -- the set of edges from $E$ with one end in $V'$ and the other in $V''$, and finally, by $G-v$ we mean the graph obtained from $G$ by removing $v$ and all its incident edges. 
Moreover, the sum of graphs $G_1=(V_1,E_1), G_2=(V_2,E_2)$ is understood as $G_1\cup G_2:=(V_1\cup V_2, E_1\cup E_2)$.

\section{Proof of Theorem~\ref{1234regTh}}

Let $G=(V,E)$ be a $d$-regular graph, $d\geq 2$, and let $I\subset V$ be an arbitrary maximal independent set in $G$. Denote $R:=V\smallsetminus I$. Let $R_1\subseteq R$ be the set of isolated vertices in $G[R]$, and set $R_2:=R\smallsetminus R_1$. Denote by $G_1,\ldots,G_p$ the components of $G[R_2]$ (each of which contains at least one edge). For every $i=1,\ldots,p$, we order the vertices of $G_i$ into a sequence $v_1,\ldots,v_{n}$ so that each $v_j$ with $j<n$ has a \emph{forward neighbour} in $G_i$, that is a neighbour $v_k$ of $v_j$ in $G_i$ with $k>j$ (this can be achieved by denoting any vertex of $G_i$ as $v_n$ and using e.g. BFS algorithm to find a spanning tree of $G_i$ rooted at $v_n$, denoting consecutive vertices encountered within the algorithm: $v_{n-1},v_{n-2},\ldots$); we denote the edge joining $v_j$ with such $v_k$ with the least index $k$ ($k>j$) the \emph{first forward edge of} $v_j$.
Analogously we define \emph{backward neighbours} of a given vertex in $G_i$.
The vertex $v_n$ shall moreover be called the \emph{last vertex} of $G_i$. 
By definition, every vertex in $R$ is incident with an edge joining it with $I$; for every $v\in R_2$ which is not the last vertex (in some component of $G[R_2]$) choose arbitrarily one such edge and denote it $e_v$ -- we shall call $e_v$ the \emph{supporting edge} of $v$.
We shall first assign initial weights $\omega(e)$ to all the edges $e$ of $G$. These shall be modified so that at the end of our construction: 
\begin{itemize}
\item[(a)] $\sigma(v)<3d$ for every $v\in R_2$;
\item[(b)] $\sigma(v)\geq 3d$ for every $v\in I$; 
\item[(b')] $\sigma(v)<4d$ for every $v\in I$ with a neighbour in $R_1$;
\item[(b'')] $\sigma(v)\in\{3d-1,4d\}$ for every $v\in R_1$,
\end{itemize}
where by $\sigma(v)$ we mean the sum at a given vertex $v$ in $G$.
Note that since $I$ and $R_1$ are stable sets and there are no edges between $R_1$ and $R_2$ in $G$, by (a), (b), (b') and (b''), potential sum conflicts shall only be possible between adjacent vertices in $R_2$ then.
We shall also require that throughout the whole construction:
\begin{itemize}
\item[(c)] 
$\omega(e)\in\{1,2,3\}$ for $e\in E(R_2)$,\\     
$\omega(e)\in\{3,4\}$ for $e\in E(I,R_2)$,\\ 
$\omega(e)\in\{2,3,4\}$ for $e\in E(I,R_1)$. 
\end{itemize}
Major concern of our weight modifyng algorithm shall be devoted to distinguishing adjacent vertices in $R_2$. 
Only in its final stage shall we adjust the sums in $R_1$ (still consistently with (a), (b), (b'), (b'') and (c)).
Initially we assign the weight:
\begin{itemize}
\item[(i)] $\omega(e)=1$, if $e$ is the first forward edge of some vertex;
\item[(ii)] $\omega(e)=2$, if $e$ is an edge of $G[R_2]$ which is not the first forward edge of any vertex;
\item[(iii)] $\omega(e)=3$, if $e$ is incident with a vertex in $I$ and is not a supporting edge;
\item[(iv)] $\omega(e)=4$, if $e$ is a supporting edge. 
\end{itemize}
Note that these weights are consistent with (c).

In the following main part of our modifying procedure we analyse and alter the sums at consecutive vertices in all the components of $G[R_2]$. Thus suppose we have already analysed all vertices in $G_1,\ldots,G_{i-1}$, and within $G_i$ -- the vertices $v_1,\ldots,v_{j-1}$ (following the rules $(1^\circ)$ -- $(3^\circ)$ specified below), hence we are about to consider the vertex $v_j$ (consistently with the vertex ordering fixed in $G_i$). While analysing this vertex: 
\begin{itemize}
\item[$(1^\circ)$] we are not allowed to modify the sums at already analysed vertices (which are fixed and shall not change till the end of the construction).
\end{itemize} 
On the other hand we wish to make some weight alterations so that: 
\begin{itemize}
\item[$(2^\circ)$] the obtained sum at $v_j$ is distinct from the sums at all the already analysed neighbours of $v_j$ in $G[R_2]$ (i.e. those in $\{v_1,\ldots,v_{j-1}\}$);
\end{itemize}
while for this aim: 
\begin{itemize}
\item[$(3^\circ)$] we are allowed to modify by $1$ the weights of the edges joining $v_j$ with its backward neighbours in $G_i$ and the weights of their supporting edges so that (c) still holds.
\end{itemize}

Before we show that we can indeed perform our modifying procedure in accordance with $(1^\circ)$ -- $(3^\circ)$,
let us observe the following. 

\begin{observation}\label{a_b_Observation} 
After analysing all vertices of $R_2$ consistently with requirements $(1^\circ)$ -- $(3^\circ)$, 
the conditions (a), (b) and (b') shall hold. 
\end{observation}

\begin{pf} 
By (iii) and (iv) all edges incident with a vertex $v\in I$ are initially weighted $3$ or $4$, 
while by $(3^\circ)$ the weight of an edge $e$ incident with $v$ can only be altered if $e$ is a supporting edge -- by $(3^\circ)$ and (c) we however must still have $\omega(e)\in\{3,4\}$ afterwards, thus (b) follows.
If $v$ has moreover a neighbour $u\in R_1$, then by (iii) we must have $\omega(uv)=3$, and this weigh 
is not modified within our procedure (cf. $(3^\circ)$), and thus (b') is fulfilled as well.

To see that (a) must also hold, note first that each edge $e$ of $G[R_2]$ can be modified at most once (consistently with $(3^\circ)$) within the algorithm, when it joins the currently analysed vertex with its backward neighbour. Therefore, for every vertex $v\in R_2$ which is not the last vertex of some component of $G[R_2]$, immediately after analysing $v$, the first forward edge of $v$ still has unchanged weight $1$ (cf. (i)). By (i) -- (iv) and $(3^\circ)$, all its remaining incident edges have in turn weights at most $3$, except for $e_v$, which has weight  $4$.  Therefore, $\sigma(v)\leq 3d-1$, and by $(1^\circ)$ this does not change till the end of  the construction. In order to prove the same holds also in the case when $v\in R_2$ is the last vertex of some component of $G[R_2]$, it is sufficient to note that then, by our construction:
$\omega(e)\leq 3$ for every edge incident with $v$, as only supporting edges can be at this point weighted $4$. 
Thus (a) follows, as by (i) and $(3^\circ)$ the edge joining $v$ with the vertex $u$ directly preceding it in the corresponding ordering cannot have weight greater than $2$ (as according to the main feature of the previously fixed orderings, this has to be the first forward edge of $u$).  
\qed
\end{pf}

Now we explain how we can perform every consecutive step of our modifying procedure, associated with a currently analysed vertex $v_j$ from component $G_i$, so that $(1^\circ)$ -- $(3^\circ)$ hold (provided that the previous steps were consistent with these rules). For this aim note first that while analysing $v_j$, the weight of every \emph{backward edge} of $v_j$ (i.e. an edge joining it with its backward neighbour in $G_i$) \textbf{can} be modified by $1$ if necessary. Indeed, suppose $e=v_kv_j$ is such an edge (i.e. $k<j$). If $e$ is not the first forward edge of $v_k$, then by (ii), $\omega(e)=2$ and by (c), $\omega(e_{v_k})\in\{3,4\}$. Thus, so that $(1^\circ)$ is obeyed, according to $(3^\circ)$, if $\omega(e_{v_k})=3$, we may change the weights of $e$ and $e_{v_k}$ to $1$ and $4$, resp., while if $\omega(e_{v_k})=4$, we may change the weights of $e$ and $e_{v_k}$ to $3$ and $3$, respectively. On the other hand, if $e$ is the first forward edge of $v_k$, then neither $\omega(e)$ nor $\omega(e_{v_k})$ have been modified thus far, hence we may modify their respective current values $1$ and $4$ to $2$ and $3$ respectively. Suppose now that $v_j$ has $b$ backward neighbours, hence also $b$ backward edges, then as each of these provides one more possible alteration of the sum at $v_j$, we altogether have $b+1$ available options for this sum (which do not influence the sums at the backward neighbours of $v_j$). Thus we may choose among these admissible alterations
such that result in $\sigma(v_j)$ distinct from sums fixed for all $b$ backward neighbours of $v_j$ in $R_2$, i.e. consistent with $(2^\circ)$.

After analysing in this manner all vertices in $R_2$ we obtain a weighting of $G$ for which, by $(1^\circ)$, $(2^\circ)$ and Observation~\ref{a_b_Observation} (which guarantees (a) and (b)), $\sigma(u)\neq \sigma(v)$ for every $uv\in E(R_2\cup I)$. 
Now we modify the sums in $R_1$ so that (b'') holds. Recall that by definition, each vertex $v\in R_1$ is only adjacent with vertices in $I$, and thus all edges  incident with such $v$ are weighted $3$ by (iii) and (c). One after another, for every $v\in R_1$ we proceed as follows.
If $\sigma(u')\geq 3d+1$ for any neighbour $u'$ of $v$ ($u'\in I$), then we change the weight of exactly one edge,
namely $u'v$ from $3$ to $2$. 
Otherwise, i.e. when due to (b) wee have $\sigma(u)=3d$ for every neighbour $u$ of $v$ ($u\in I$), we change the weight of $uv$ from $3$ to $4$ for all $u\in N(v)\subseteq I$. Note that in the both cases none of (a), (b) and (b') shall be violated, while we shall attain: 
 $\sigma(v)\in\{3d-1,4d\}$. After processing in this manner consecutively all vertices in $R_1$, all neighbours in $G$ shall finally be sum-distinguished, as vertices in $R_1$ are only adjacent with those in $I$, cf. (b), (b') and (b'').
\qed

\section{Proof of Theorem~\ref{123largeregTh}\label{SectionProofLargeD}}

\subsection{Tools}

The proof of Theorem~\ref{123largeregTh} relies heavily on a random distribution of vertices and edges of a given graph 
to subsets with carefully predefined proportions. For this aim we shall however also make use of Corollary~\ref{QuarterDecompositionLemma} below implied by the following straightforward deterministic observation from~\cite{PrzybyloStandard22}, and possibly many other sources.
\begin{observation}
\label{EvenDecomposition}
Every graph $G=(V,E)$ can be edge decomposed into two subgraphs $G_1, G_2$ so that for each $v\in V$ and $i\in{1,2}$:
\begin{equation}\label{EQ_EulerianDecomposition}
d_{G_i}(v)\in \left[\frac{d_G(v)}{2}-1,\frac{d_G(v)}{2}+1\right].
\end{equation} 
\end{observation}

\begin{corollary}\label{QuarterDecompositionLemma}
Every  graph $G_1=(V_1,E_1)$ has a subgraph $G'_1$ such that for each  $v\in V_1$:
\begin{equation}\label{dG'1}
d_{G'_1}(v)\in\left[\frac{9}{16}d_{G_1}(v)-3,\frac{9}{16}d_{G_1}(v)+3\right].
\end{equation}
\end{corollary}

\begin{pf}
Such a subgraph can be constructed via four times repeated application of Observation~\ref{EvenDecomposition},
first to $G_1$ to obtain say $G^{(1)}_2$ and $G^{(2)}_2$, then to $G^{(2)}_2$ to obtain $G^{(1)}_3$ and $G^{(2)}_3$, next to $G^{(2)}_3$ to obtain $G^{(1)}_4$ and $G^{(2)}_4$, and finally to $G^{(2)}_4$ to get $G^{(1)}_5$ and $G^{(2)}_5$. It is then straightforward to verify that~(\ref{EQ_EulerianDecomposition}) implies that $G'_1:=G^{(1)}_2\cup G^{(2)}_5$ complies with our requirements. 
\qed
\end{pf}
For random arguments we shall mostly use the symmetric variant of the Lov\'asz Local Lemma, see e.g.~\cite{AlonSpencer} and the Chernoff Bound, see e.g.~\cite{JansonLuczakRucinski}. 
\begin{theorem}[\textbf{The Local Lemma}]
\label{LLL-symmetric}
Let $A_1,A_2,\ldots,A_n$ be events in an arbitrary pro\-ba\-bi\-li\-ty space.
Suppose that each event $A_i$ is mutually independent of a set of all the other
events $A_j$ but at most $D$, and that $\mathbf{Pr}(A_i)\leq p$ for all $1\leq i \leq n$. If
$$ p \leq \frac{1}{e(D+1)},$$
then $ \mathbf{Pr}\left(\bigcap_{i=1}^n\overline{A_i}\right)>0$.
\end{theorem}
\begin{theorem}[\textbf{Chernoff Bound}]\label{ChernofBoundTh}
For any $0\leq t\leq np$,
$$\mathbf{Pr}\left(\left|{\rm BIN}(n,p)-np\right|>t\right)<2e^{-\frac{t^2}{3np}}$$
where ${\rm BIN}(n,p)$ is the sum of $n$ independent Bernoulli variables, each equal to $1$ with probability $p$ and $0$ otherwise.
\end{theorem}

Finally, the following technical observation shall be useful repeatedly throughout the proof of Theorem~\ref{123largeregTh} while applying the local lemma. 

\begin{observation}
For every $x\geq 10^8$,
\begin{eqnarray}
&&2e^{-\frac{x}{2.45\cdot 10^6}} <\frac{1}{2ex^2}; \label{TechIneq1}\\
&&2e^{-\frac{x}{4.9\cdot 10^6}} < \frac{1}{ex}. \label{TechIneq2}
\end{eqnarray}
\end{observation}

\begin{pf}
Note first that (\ref{TechIneq1}) is directly implied by inequality (\ref{TechIneq2}) -- it is sufficient to square the both sides of (\ref{TechIneq2}).
To prove inequality~(\ref{TechIneq2})  it is however equivalently sufficient to show that $f(x)>0$ for $x\geq 10^8$ where
$$f(x):=\frac{x}{4.9\cdot 10^6}-\ln(2ex).$$
This in turn holds since $f'(x) = \frac{1}{4.9\cdot 10^6}-\frac{1}{x}>0$ for $x > 4.9\cdot 10^6$ and 
$f(10^8) = \frac{100}{4.9}-\ln\left(2e10^8\right) = \frac{20\cdot4.9+2}{4.9}-\ln2-1-8\ln10 
> 20.4-0.7-1-8\cdot 2.31> 0$.
\qed
\end{pf}

\subsection{Random vertex and edge partitions}

Let $G=(V,E)$ be a $d$-regular graph with $d\geq 10^8$.
\begin{claim}\label{ClaimV0}
We can choose a subset $V_0\subseteq V$ such that for every $v\in V$:
\begin{equation}\label{dV0}
\left|d_{V_0}(v)-0.05d\right|\leq 3\cdot10^{-4}d.
\end{equation}
\end{claim}

\begin{pf}
Randomly and independently we place every vertex from $V$ in $V_0$ with probability $0.05$. 
Denote by $A_{1,v}$ the event that (\ref{dV0}) does not hold for a given $v\in V$.
By the Chernoff Bound, i.e. Lemma~\ref{ChernofBoundTh}, and inequality~(\ref{TechIneq1}):
$$\mathbf{Pr}\left(A_{1,v}\right) < 2e^{-\frac{(3\cdot10^{-4}d)^2}{3\cdot0.05d}} 
= 2e^{-\frac{d}{\frac53\cdot 10^6}}  <\frac{1}{ed^2}.$$ 
As every event $A_{1,v}$ is mutually independent of all other events $A_{1,u}$ except those where $u$ share a common neighbour with $v$, i.e. all except at most $d(d-1)<d^2-1$ events, by the Lov\'asz Local Lemma, i.e. Lemma~\ref{LLL-symmetric}, with positive probability none of the events $A_{1,v}$ holds, and thus $V_0$ as desired must exist.
\qed
\end{pf}

Fix any $V_0$ consistent with Claim~\ref{ClaimV0}, and denote:
$$V_1:=V\smallsetminus V_0,~~~~
G_1:=G[V_1] {\rm~~~~and~~~~} G_0:=G[V_0],$$  
hence by~(\ref{dV0}), for every $v\in V$:
\begin{equation}\label{dV1}
\left|d_{V_1}(v) - 0.95d\right|\leq 3\cdot10^{-4}d.
\end{equation}

We shall first fix sums for all vertices in $V_1$ (containing great majority of all the vertices), keeping these relatively small, and using weights of some of the edges between $V_1$ and $V_0$ for some necessary adjustments.
By Corollary~\ref{QuarterDecompositionLemma} there is a subgraph $G'_1$ of $G_1$ such that~(\ref{dG'1}) holds for  every $v\in V_1$.
By~(\ref{dG'1}) and~(\ref{dV1}) we thus obtain that for every $v\in V_1$:
\begin{equation}\label{dG'1_2}
0.5d < \frac{9}{16}(0.95-0.0003)d-3\leq d_{G'_1}(v) \leq \frac{9}{16}(0.95+0.0003)d+3 < 0.54 d. 
\end{equation}

Let
 $$c_1:V_1\to\{1,2,\ldots,10^4\}$$ 
be an auxiliary assignment of integers to the vertices of $G_1$. Denote:
\begin{eqnarray}
E': &=& \left\{uv\in E(G'_1): c_1(u)+c_1(v)\geq 10^4+2\right\}, \label{DefinitionOfE'}\\
E'': &=& \left\{uv\in E(G'_1): c_1(u)+c_1(v)\leq 10^4+1\right\} \label{DefinitionOfE''}
\end{eqnarray}
(the edges in $E'$ shall be the only ones in $G_1$ with weight $3$ assigned, while the remaining ones shall be weighted $1$ -- this, combined with Claim 2 below shall assure convenient sums distribution in neighbourhoods of all vertices in $V_1$), and for each $i=1,2,\ldots,10^4$:
$$V_{1,i}: = \left\{v\in V_1: c_1(v)=i\right\}.$$

\begin{claim}\label{ClaimV1iE'}
We may choose $c_1$ so that for every $i\in \{1,2,\ldots,10^4\}$ and each $v\in V_{1,i}$:
\begin{eqnarray}
&& \left| d_{V_{1,i}}(v) - 10^{-4}d_{V_1}(v) \right| \leq 11\cdot 10^{-6}d; \label{dV1i}\\
&& \left| d_{E'}(v) - (i-1)10^{-4}d_{G'_1}(v) \right| \leq 6\cdot 10^{-4}d. \label{dE'}
\end{eqnarray}
\end{claim}

\begin{pf}
We choose $c_1:V_1\to\{1,2,\ldots,10^4\}$ randomly by independently assigning every vertex $v\in V_1$ its value $c_1(v)$ from the set $\{1,2,\ldots,10^4\}$, each with equal probability. 
For every $v\in V_1$ we denote 
by  $A_{2,v}$ and $A_{3,v}$ the events that (\ref{dV1i}) and that (\ref{dE'}) does not hold, respectively.
Let $v\in V_1$. As by~(\ref{dV1}) we have $11\cdot 10^{-6}d\leq 10^{-4}d_{V_1}(v)$, by the Chernoff Bound, (\ref{dV1}) and (\ref{TechIneq1}) we obtain that: 
\begin{eqnarray}
\mathbf{Pr}\left(A_{2,v}\right)<2e^{-\frac{(11\cdot 10^{-6}d)^2}{3\cdot 10^{-4}d_{V_1}(v)}} \leq 2e^{-\frac{121\cdot 10^{-12}d^2}{3\cdot 10^{-4}(0.95+3\cdot10^{-4})d}} = 
2e^{-\frac{d}{\frac{285.09}{121}\cdot 10^{6}}} 
< \frac{1}{2ed^2}. \label{PrA2v}
\end{eqnarray}

Note further that for every $i\geq 2$, $i\leq 10^4$, by~(\ref{DefinitionOfE'}) and~(\ref{DefinitionOfE''}):
\begin{eqnarray}
&&\mathbf{Pr}\left(A_{3,v}~|~c_1(v)=i\right) \nonumber\\
&=& \mathbf{Pr}\left(\left|{\rm BIN}\left(d_{G'_1}(v),(i-1)10^{-4}\right) - (i-1)10^{-4}d_{G'_1}(v)\right|>6\cdot 10^{-4}d\right)   \nonumber\\
&=&\mathbf{Pr}\left(\left|d_{E'}(v)-(i-1)10^{-4}d_{G'_1}(v)\right| > 6\cdot 10^{-4}d ~|~ c_1(v)=i\right)  \nonumber\\
&=&\mathbf{Pr}\left(\left|d_{G'_1}(v)-d_{E''}(v)-(i-1)10^{-4}d_{G'_1}(v)\right| > 6\cdot 10^{-4}d ~|~ c_1(v)=i\right)  \nonumber\\
&=& \mathbf{Pr}\left(\left|(10^4-i+1)10^{-4}d_{G'_1}(v)-d_{E''}(v)\right| > 6\cdot 10^{-4}d ~|~ c_1(v)=i\right)  \nonumber\\
&=& \mathbf{Pr}\left(\left|{\rm BIN}\left(d_{G'_1}(v),(10^4-i+1)10^{-4}\right) - 
(10^4-i+1)10^{-4}d_{G'_1}(v)\right|>6\cdot 10^{-4}d\right)  \nonumber\\
&=& \mathbf{Pr}\left(A_{3,v}~|~c_1(v)=10^4-i+2\right). \label{A3vSymmetry}
\end{eqnarray}
Now for every fixed $13\leq i\leq 0.5\cdot 10^4+1$ 
(as then by~(\ref{dG'1_2}), 
$(i-1)10^{-4}d_{G'_1}(v) > 12\cdot 10^{-4} \cdot 0.5d = 6\cdot 10^{-4}d$), by the Chernoff Bound, (\ref{dG'1_2}) and (\ref{TechIneq1}) we obtain: 
\begin{equation}
\mathbf{Pr}\left(A_{3,v}~|~c_1(v)=i\right) < 2e^{-\frac{(6\cdot 10^{-4}d)^2}{3\cdot (i-1)10^{-4}d_{G'_1}(v)}} 
<  2e^{-\frac{(6\cdot 10^{-4}d)^2}{3\cdot 0.5 \cdot 0.54d}} 
=  2e^{-\frac{d}{2.25 \cdot 10^6}} 
<  \frac{1}{2ed^2}. \label{A3vMostCases}
\end{equation}
For $i=1$, by the definition of $E'$ and $c_1$, we trivially have:  
\begin{equation}\label{A3vCase0}
\mathbf{Pr}\left(A_{3,v}~|~c_1(v)=1\right) = 0.
\end{equation}
For $i\in\{2,3,\ldots,12\}$ in turn (as then by~(\ref{dG'1_2}), $(i-1)10^{-4}d_{G'_1}(v) > 0.5\cdot 10^{-4}d$), by the Chernoff Bound and (\ref{TechIneq1}):
\begin{eqnarray}
 \mathbf{Pr}\left(A_{3,v}~|~c_1(v)=i\right) 
 &\leq& \mathbf{Pr}\left(\left|d_{E'}(v)-(i-1)10^{-4}d_{G'_1}(v)\right| > 0.5\cdot 10^{-4}d ~|~ c_1(v)=i\right) \nonumber\\ 
& <& 2e^{-\frac{( 0.5\cdot 10^{-4}d)^2}{3\cdot (i-1)10^{-4}d_{G'_1}(v)}} 
< 2e^{-\frac{( 0.5\cdot 10^{-4}d)^2}{3\cdot 11\cdot 10^{-4}d}} 
= 2e^{-\frac{d}{\frac{33}{25}\cdot 10^{6}}}  
<  \frac{1}{2ed^2}. \label{A3vSmallCases}
 \end{eqnarray}
By (\ref{A3vMostCases}), (\ref{A3vSmallCases}), (\ref{A3vSymmetry}), (\ref{A3vCase0})  and the law of total probability,
\begin{equation}\label{PrA3v}
\mathbf{Pr}\left(A_{3,v}\right) < \frac{1}{2ed^2}.
\end{equation}
Let $\Delta_1$ be the maximum degree of $G_1$.
As every event $A_{2,v}$ and every event $A_{3,v}$ is mutually independent of all other events $A_{2,u}$ and $A_{3,u}$ except possibly those where $u$ is at distance at most $2$ from $v$ in $G_1$, i.e. all except at most $2\Delta_1^2+1<2d^2-1$, by~(\ref{PrA2v}), (\ref{PrA3v}) and the Lov\'asz Local Lemma, with positive probability none of the events $A_{2,v}$, $A_{3,v}$ holds, and thus $c_1$ as desired must exist.
\qed
\end{pf}

\begin{claim}\label{ClaimE1inV0andV1}
We may choose  a set of edges $E_1\subseteq E(V_1,V_0)$ such that:
\begin{eqnarray}
&& \left| d_{E_1}(u) - 0.08d_{V_0}(u) \right| \leq 5\cdot10^{-5}d {\rm ~~~~for~~~~} u\in V_1; \label{dE1inV1}\\
&& \left| d_{E_1}(v) - 0.08d_{V_1}(v) \right| \leq 10^{-3}d  {\rm ~~~~for~~~~} v\in V_0. \label{dE1inV0}
\end{eqnarray}
\end{claim}

\begin{pf}
Suppose we randomly and independently place every edge from $E(V_0,V_1)$ in $E_1$ with probability $0.08$.
For every $u\in V_1$ and $v\in V_0$, denote 
by $A_{4,u}$ and $A_{5,v}$ the events that (\ref{dE1inV1}) and that (\ref{dE1inV0}) does not hold, respectively.
Then by the Chernoff Bound, (\ref{dV0}), (\ref{dV1}) and (\ref{TechIneq2}) we obtain that for every $u\in V_1$:
\begin{equation}\label{PrA4v}
\mathbf{Pr}\left(A_{4,u}\right) < 2e^{-\frac{(5\cdot 10^{-5}d)^2}{3\cdot 0.08d_{V_0}(u)}} 
\leq  2e^{-\frac{(5\cdot 10^{-5}d)^2}{3\cdot0.08\cdot0.0503d}} 
=  2e^{-\frac{d}{4.8288 \cdot 10^{6}}} 
< \frac{1}{ed},
\end{equation}
while for each $v\in V_0$:
\begin{equation}\label{PrA5v}
\mathbf{Pr}\left(A_{5,v}\right) < 2e^{-\frac{(10^{-3}d)^2}{3\cdot0.08d_{V_1}(v)}} \leq 
 2e^{-\frac{(10^{-3}d)^2}{3\cdot0.08\cdot0.9503d}} = 
 2e^{-\frac{d}{24\cdot 9503}} <
\frac{1}{ed}.
\end{equation}
As every event $A_{4,v}$ and every event $A_{5,v}$ is mutually independent of all other events $A_{4,u}$ and $A_{5,u}$ except possibly those where $u$ is at distance at most $1$ from $v$ in the graph induced by the edges of $E(V_0,V_1)$, i.e. (by (\ref{dV0}) and (\ref{dV1})) all except at most $\max_{v\in V} d_{E(V_0,V_1)}(v) < d-1$, by~(\ref{PrA4v}), (\ref{PrA5v}) and the Lov\'asz Local Lemma, with positive probability none of the events $A_{4,v}$, $A_{5,v}$ holds, and thus $E_1$ as desired must exist.
\qed
\end{pf}

Carefully designed weight adjustments of the edges in $E_1$ shall be used to supplement roughly even distribution of sums assured by Claim~\ref{ClaimV1iE'}, and consequently to provide sum distinction of the neighbours in $V_1$. 
Out of the remaining edges in $E(V_0,V_1)$ we shall next choose a set $E_0$ of edges with special features, which 
shall partition $V_0$ into 5 subsets (with decreasing sums, all however larger than the ones in $V_1$), and shall facilitate the use of a modification of Kalkowski's algorithm in $G_0$ to distinguish the remaining neighbours in $G$.
Denote:
\begin{equation}\label{E*definition}
E^*:= E(V_0,V_1)\smallsetminus E_1.
\end{equation}
Consider an assignment:
 $$c_0:V_0\to\{0,1,2,3,4\}$$ 
and the following partition of $V_0$ induced by it: 
$$V_{0,j}:=\left\{v\in V_0: c_0(v)=j\right\}, ~~~~j=0,1,2,3,4.$$

\begin{claim}\label{ClaimE0inV0andV1anddV0i}
We may choose $E_0\subseteq E^*$ and $c_0:V_0 \to \{0,1,2,3,4\}$ such that:
\begin{eqnarray}
&& \left| d_{E_0}(v) - 0.2c_0(v)d_{E^*}(v) \right| \leq 10^{-3}d {\rm ~~~~for~~~~} v\in V_0; \label{dE0inV0}\\
&& \left| d_{V_{0,c_0(v)}}(v)- 0.2d_{V_0}(v)\right| \leq 10^{-3}d {\rm ~~~~for~~~~} v\in V_0; \label{dV0i}\\
&& \left| d_{E_0}(u) - 0.4d_{E^*}(u) \right| \leq 2\cdot10^{-4}d  {\rm ~~~~for~~~~} u\in V_1. \label{dE0inV1}
\end{eqnarray}
\end{claim}

\begin{pf} First for every $v\in V_0$ choose randomly and independently an integer in $\{0,1,2,3,4\}$, each with equal probability, and denote it by $c_0(v)$. Then include every edge $uv\in E^*$ with $v\in V_0$ in $E_0$ randomly and independently with probability $0.2c_0(v)$. For every $v\in V_0$ and $u\in V_1$, denote by
$A_{6,v}$, $A_{7,v}$ and $A_{8,u}$ the events that (\ref{dE0inV0}), (\ref{dV0i}) and (\ref{dE0inV1}) does not hold, respectively.
Let $v\in V_0$. 
First note that by the Chernoff Bound, (\ref{dE1inV0}), (\ref{dV1}) and (\ref{TechIneq1}), for every fixed $j\in\{1,2,3,4\}$:
\begin{eqnarray}
\mathbf{Pr}\left(A_{6,v} ~|~ c_0(v)=j\right) &<& 2e^{-\frac{(10^{-3}d)^2}{3\cdot 0.2j\cdot d_{E^*}(v)}} \leq
2e^{-\frac{(10^{-3}d)^2}{0.6 j (0.92d_{V_1}(v)+ 10^{-3}d)}} \nonumber\\ 
&\leq& 2e^{-\frac{(10^{-3}d)^2}{2.4\cdot (0.92\cdot 0.9503d+ 10^{-3}d)}}  
= 2e^{-\frac{d}{2100662.4}}  < \frac{1}{2ed^2}. \label{PrA6v1234}
\end{eqnarray}
Moreover, if $c_0(v)=0$, then we trivially have: $d_{E_0}(u)=0$, and hence:  
\begin{equation}
\mathbf{Pr}\left(A_{6,v} ~|~ c_0(v)=0\right) = 0. \label{PrA6v0}
\end{equation}
By (\ref{PrA6v1234}), (\ref{PrA6v0}) and the law of total probability we thus obtain that:
\begin{equation}
\mathbf{Pr}\left(A_{6,v}\right) < \frac{1}{2ed^2}. \label{PrA6v}
\end{equation}

Further,  by the Chernoff Bound, (\ref{dV0}) and (\ref{TechIneq1}), 
\begin{eqnarray}
\mathbf{Pr}\left(A_{7,v} \right) &<& 2e^{-\frac{(10^{-3}d)^2}{3\cdot 0.2d_{V_0}(v)}} \leq
2e^{-\frac{(10^{-3}d)^2}{ 0.6\cdot 0.0503 d}} =
2e^{-\frac{ d}{30180}}  <
\frac{1}{2ed^2}.  \label{PrA7v}
\end{eqnarray}

Let now $v\in V_1$.
Note that as our choices are independent, then every edge from $E^*$ which is incident with $v$ is in fact independently chosen to $E_0$ with probability $$\frac{1}{5}\cdot(0+0.2+0.4+0.6+0.8)=0.4,$$ and hence, by the Chernoff Bound,
(\ref{dE1inV1}), (\ref{dV0}) and (\ref{TechIneq1}), 
\begin{eqnarray}
\mathbf{Pr}\left(A_{8,v}\right) &<& 2e^{-\frac{(2\cdot10^{-4}d)^2}{3\cdot 0.4d_{E^*}(v)}} \leq
2e^{-\frac{(2\cdot10^{-4}d)^2}{1.2\cdot (0.92d_{V_0}(v)+5\cdot 10^{-5}d)}} \leq 
2e^{-\frac{(2\cdot10^{-4}d)^2}{1.2\cdot (0.92\cdot 0.0503d+5\cdot 10^{-5}d)}} \nonumber\\ 
& = &2e^{-\frac{d}{1389780}}  < \frac{1}{2ed^2}. \label{PrA8v}
\end{eqnarray}

It is easy to notice that each of the events $A_{6,v}$, $A_{7,v}$, $A_{8,v}$ is mutually independent of all but (much) less than $2d^2$ other events of such types, and thus by (\ref{PrA6v}), (\ref{PrA7v}), (\ref{PrA8v}) and the Lov\'asz Local Lemma, with positive probability none of the events $A_{6,v}$, $A_{7,v}$, $A_{8,v}$ holds, and hence there must exist $c_0$ and $E_0$ as required.
\qed
\end{pf}

\subsection{Initial weighting}

We define an initial edge 3-weighting $\omega_0$ of $G$ as follows:
$$\omega_0(e):=\left\{ \begin{array}{lll} 
1,&{\rm if} &e\in E(V_1)\smallsetminus E';\\
2,&{\rm if}& e\in E_1\cup E_0\cup E(V_0);\\
3,&{\rm if}& e\in E'\cup E(V_0,V_1)\smallsetminus(E_0\cup E_1).
\end{array}\right.$$
We shall further modify only the weights of the edges in $E_1$ (increasing some of these to $3$, to adjust the sums in $V_1$) and of the edges in $E(V_0)$ (possibly changing some of them by $1$ in order to distinguish sums within $V_0$ in the final part of our construction).

Note that at this point for every $v\in V_1$:
\begin{eqnarray}
\sigma(v) 
&=& 3\cdot d_{E'}(v)+1\cdot \left(d_{V_1}(v)-d_{E'}(v)\right)+2\cdot \left(d_{E_0}(v)+d_{E_1}(v)\right)+3\cdot \left(d_{V_0}(v)-d_{E_0}(v)-d_{E_1}(v)\right)\nonumber\\
&=& d+2d_{E'}(v)+1\cdot \left(d_{E_0}(v)+d_{E_1}(v)\right)+2\cdot \left(d_{V_0}(v)-d_{E_0}(v)-d_{E_1}(v)\right) \nonumber\\
&=& d+2d_{E'}(v)+2d_{V_0}(v)-d_{E_0}(v)-d_{E_1}(v), \nonumber
\end{eqnarray}
and hence, if $v\in V_{1,i}$ for some $i\in\{1,2,\ldots,10^4\}$, then by (\ref{dE'}), (\ref{dE0inV1}), (\ref{E*definition}), (\ref{dG'1}), (\ref{dE1inV1}), (\ref{dV0}):
\begin{eqnarray}
&&\sigma(v) \in\nonumber\\
&&\Big[d+2\left((i-1)10^{-4}d_{G'_1}(v)-6\cdot 10^{-4}d\right)+2d_{V_0}(v)
-\left(0.4\left(d_{V_0}(v)-d_{E_1}(v)\right)+2\cdot10^{-4}d\right)-d_{E_1}(v), \nonumber\\
&&d+2\left((i-1)10^{-4}d_{G'_1}(v)+6\cdot 10^{-4}d\right)+2d_{V_0}(v)
-\left(0.4\left(d_{V_0}(v)-d_{E_1}(v)\right)-2\cdot10^{-4}d\right)-d_{E_1}(v)\Big]\nonumber\\
&\subseteq& \bigg[d+2(i-1)10^{-4}\left(\frac{9}{16}(d-d_{V_0}(v))-3\right)+1.6d_{V_0}(v)-0.6d_{E_1}(v) -0.0014d, \nonumber\\
&&d+2(i-1)10^{-4}\left(\frac{9}{16}(d-d_{V_0}(v))+3\right)+1.6d_{V_0}(v)-0.6d_{E_1}(v) +0.0014d\bigg] \nonumber\\
&\subseteq& \bigg[d+2(i-1)10^{-4}\left(\frac{9}{16}(d-d_{V_0}(v))-3\right)+1.6d_{V_0}(v) 
-0.6\left(0.08d_{V_0}(v)+5\cdot10^{-5}d\right) -0.0014d, \nonumber\\
&&d+2(i-1)10^{-4}\left(\frac{9}{16}(d-d_{V_0}(v))+3\right)+1.6d_{V_0}(v)
-0.6\left(0.08d_{V_0}(v)-5\cdot10^{-5}d\right) +0.0014d\bigg] \nonumber\\
&=& \bigg[\left(1+\frac{9}{8}(i-1)10^{-4}\right)d+\left(1.552-\frac{9}{8}(i-1)10^{-4}\right)d_{V_0}(v) -0.00143d - 6(i-1)10^{-4}, \nonumber\\
&&\left(1+\frac{9}{8}(i-1)10^{-4}\right)d+\left(1.552-\frac{9}{8}(i-1)10^{-4}\right)d_{V_0}(v) + 0.00143d + 6(i-1)10^{-4}\bigg]\nonumber\\
&\subseteq& \bigg[\left(1+\frac{9}{8}(i-1)10^{-4}\right)d+\left(1.552-\frac{9}{8}(i-1)10^{-4}\right)\left(0.05d-0.0003d\right) 
-0.00143d - 6(i-1)10^{-4}, \nonumber\\
&& \left(1+\frac{9}{8}(i-1)10^{-4}\right)d+\left(1.552-\frac{9}{8}(i-1)10^{-4}\right)\left(0.05d+0.0003d\right) 
+ 0.00143d + 6(i-1)10^{-4}\bigg] \nonumber\\
&=& \big[\left(1.0776+1.06875(i-1)10^{-4}\right)d - 0.0018956d 
+ 0.0003375(i-1)10^{-4}d - 6(i-1)10^{-4}, \nonumber \\ 
&& \left(1.0776+1.06875(i-1)10^{-4}\right)d + 0.0018956d 
- 0.0003375(i-1)10^{-4}d + 6(i-1)10^{-4}\big] \nonumber\\
&\subseteq& \big[\left(1.0776+1.06875(i-1)10^{-4}\right)d - 0.0018956d, 
\left(1.0776+1.06875(i-1)10^{-4}\right)d + 0.0018956d\big]. \nonumber\\\label{FirstSv}
\end{eqnarray}

Let  
$\Delta_2:=\max_{1\leq i \leq 10^4} \Delta(G[V_{1,i}])$.
By~(\ref{dV1i}) and~(\ref{dV1}), 
\begin{eqnarray}
\Delta_2 &\leq& 10^{-4}\cdot 0.9503d+11\cdot 10^{-6}d 
= 10^{-4}\cdot 1.0603 d. \label{Delta2}
\end{eqnarray}

For every $i=1,2,\ldots,10^4$ we arbitrarily choose a proper vertex colouring of $G[V_{1,i}]$: 
\begin{equation}\label{DefinitionOfc1i}
c_{1,i}: V_{1,i}\to\{0,1,\ldots,\Delta_2\}.
\end{equation}

Now we modify weights of some of the edges in $E_1$ by adding $1$ to them (hence switching their weights from $2$ to $3$) so that for every $i\in \{1,2,\ldots,10^4\}$ and each $v\in V_{1,i}$: 
\begin{equation}
\sigma(v)=\left\lfloor \left(1.0776+1.06875(i-1)10^{-4}\right)d + 0.0018956d\right\rfloor+c_{1,i}(v). \label{SecondSv}
\end{equation} 
This is feasible by (\ref{FirstSv}), as by~(\ref{dE1inV1}), (\ref{dV0}) and~(\ref{Delta2}):
\begin{eqnarray}
d_{E_1}(v) &\geq& 0.08d_{V_0}(v) -5\cdot 10^{-5}d \geq 0.08 \cdot \left(0.05d - 3\cdot 10^{-4}d\right) 
-5\cdot 10^{-5}d \nonumber\\ 
&=& 0.003926 d> 0.00389723 d \geq 2\cdot 0.0018956d + \Delta_2. \nonumber
\end{eqnarray}
We denote the obtained weighting of the edges of $G$ by $\omega_1$. As a result, by (\ref{SecondSv}) and the definitions of $c_{1,i}$, neighbours are sum-distinguished within every $V_{1,i}$ in $G$, i.e. for each $i=1,2,\ldots,10^4$ 
and every edge $uv\in E(V_{1,i})$ we have $\sigma(u)\neq \sigma(v)$. In fact however, all neighbours in $V_1$ are at this point sum-distinguished, as no conflicts are possible between distinct sets $V_{1,i}$. To see this it is sufficient to observe that for each $1\leq i < 10^4$ and any $u\in V_{1,i}$ and $v\in V_{1,i+1}$, by  (\ref{SecondSv})  and~(\ref{Delta2}) we now have: 
\begin{eqnarray}
\sigma(u) &\leq& \left(1.0776+1.06875(i-1)10^{-4}\right)d + 0.0018956d + 10^{-4}\cdot 1.0603 d \nonumber\\
&<& \left(1.0776+1.06875(i-1)10^{-4}\right)d + 0.0018956d + 1.06875\cdot 10^{-4}d - 1 \nonumber\\ 
&<&  \left\lfloor\left(1.0776+1.06875\cdot i\cdot 10^{-4}\right)d + 0.0018956d\right\rfloor\leq \sigma(v). \label{V1iSmallerV1i+1}
\end{eqnarray}

\subsection{Setting sums in $V_0$}

We shall now modify the sums in $V_0$ by altering weights of some of the edges in $E(V_0)$ so that there are no sum conflicts within the sets $V_{0,i}$. 
First however we choose for every vertex $v\in V_0$ a set $E_v\subset E(V_0)$ of its \emph{personal} incident edges, whose weights' modifications shall settle the final sum at $v$ (up to an additive factor of $1$). These sets shall satisfy the following two features:
\begin{eqnarray}
E_u\cap E_v = \emptyset &{\rm ~~for~~}& u,v\in V_0, u\neq v; \label{emptyintersectionFeature}\\ 
|E_v| \geq 0.5d_{V_0}(v)-1 &{\rm ~~for~~}& v\in V_0. \label{halfFeature}
\end{eqnarray} 
We define these for each component of $G_0$ separately.
Suppose $H$ is any such component. We add a new vertex $u$ (if necessary) and join it with single edges with all vertices of odd degrees in $H$, thus obtaining an Eulerian graph $F$ of $H$. Then we traverse the edges of $F$ along any its Eulerian tour, starting at $u$ (or at any other vertex if $H$ was itself Eulerian), and temporarily direct these edges
 consistently with our direction of movement along the Eulerian tour. We then remove the vertex $u$ (if we priory had to add it) and for every vertex $v$ we define $E_v$ as the set of edges in $H$ outgoing from $v$. It is straightforward to verify then that such sets $E_v$ meet our requirements~(\ref{emptyintersectionFeature}) and~(\ref{halfFeature}).

We now arbitrarily arrange the vertices of $V_0$ into a sequence $v_1,v_2,\ldots, v_k$
and analyse them one after another. Once we reach a given vertex $v$ we associate to it a set 
$$S_v\in \mathbb{S}:=\left\{\left\{2i,2i+1\right\}:i\in\mathbb{Z}\right\}$$ 
distinct from all sets $S_u$ already assigned to neighbours $u$ of $v$ from $V_{0,c_0(v)}$ (hence also disjoint from 
them due to the definition of $\mathbb{S}$) and guarantee that ever since this moment the sum at $v$ belongs to this set 
(note this shall guarantee that neighbours within every $V_{0,j}$ shall be sum-distinguished for $j=0,1,2,3,4$ at the end of our algorithm).  To achieve such a goal for the currently analysed vertex $v$, we admit, if necessary, one of the following two modifications of the weight of every edge $uv\in E_v$: 
\begin{itemize}
\item increasing the weight of $uv$ to $3$ (from $2$) if $u$ was not yet analysed, i.e. $v$ precedes $u$ in the chosen ordering; 
\item otherwise, changing the weight of $uv$ to $1$ or $3$ so that as a result we still have $\sigma(u)\in S_u$. 
\end{itemize}
Note that exactly one of the two options of changing the weight of $uv$ (either to $1$ or to $3$) from the second point 
above is always available. Note moreover that in order to achieve our goal we cannot set the sum at $v$ to be equal to 
a number from any two-element set $S_u$ already associated to a neighbour $u$ of $v$ from $V_{0,c_0(v)}$, i.e. we must avoid at most $2d_{V_{0,c_0(v)}}(v)$ integers, while the admitted modifications of weights of the edges in $E_v$ yield
\begin{eqnarray}
|E_v|+1 \geq 0.5d_{V_0}(v) &=& 2\left(0.2 d_{V_0}(v) + 10^{-3}d\right) + \left(0.1 d_{V_0}(v) - 2\cdot 10^{-3}d\right) \nonumber\\
&\geq&2d_{V_{0,c_0(v)}}(v) +\left(0.1\cdot \left(0.05d-3\cdot 10^{-4}d\right) - 2\cdot 10^{-3}d\right)  
> 
 2d_{V_{0,c_0(v)}}(v) \nonumber
\end{eqnarray}
potential options for the sum at $v$ (cf. (\ref{halfFeature}), (\ref{dV0i}) and (\ref{dV0})).
We then choose any of these options, say $s^*$ which does not belong to any (already fixed) $S_u$ with 
$u\in N_{V_{0,c_0(v)}}(v)$ and which requires modification of at most $2d_{V_{0,c_0(v)}}(v)$ weights of the edges incident with $v$. 
We then perform these at most $2d_{V_{0,c_0(v)}}(v)$ admitted weight modifications in $E_v$ so that $\sigma(v)=s^*$ 
(note 
the weights modified in this step are the only ones incident with $v$ with a possible value $1$) and choose as $S_v$ the only element of $\mathbb{S}$ containing $s^*$.
After analysing all vertices in $V_0$ we obtain our final edge $3$-weighting $\omega_2$ of $G$, such that there are no sum conflicts between neighbours within any $V_{0,i}$, $i=0,1,2,3,4$. Moreover, as within the algorithm above, for every $v\in V_0$ we might have had at most $2d_{V_{0,c_0(v)}}(v)$ edges weighted $1$ immediately after fixing $S_v$ and assuring that the sum at $v$ belongs to this set, and this sum at $v$ could change only by $1$ throughout the rest of the algorithm  (so that $\sigma(v)\in S_v$), hence at the end of its execution we still have for every $v\in V_{0,i}$: 
\begin{equation}\label{FirstSvEstimation}
\sigma(v)\geq 1 \cdot 2d_{V_{0,i}}(v) + 2 \cdot\left(d_{V_0}(v)-2d_{V_{0,i}}(v)\right) -1 + \sum_{u\in N_{V_1}(v)}\omega_2(uv).
\end{equation}

\subsection{Final calculations}

It remains to show that there are not possible any sum conflicts between neighbours from different sets $V_{0,i}$ and between neighbours from $V_0$ and $V_1$.
Note first that for every $i=0,1,2,3,4$ and each $v\in V_{0,i}$,  by (\ref{FirstSvEstimation}), (\ref{dE0inV0}), (\ref{dV0i}), (\ref{dE1inV0}) and (\ref{dV0}): 
\begin{eqnarray}
\sigma(v)
&\geq& 2 d_{V_0}(v)-2d_{V_{0,i}}(v) -1+ 2 \left(d_{E_1}(v)+d_{E_0}(v)\right)+ 3\left(d_{V_1}(v) - d_{E_1}(v) - d_{E_0}(v)\right) \nonumber\\
&=& 2d+d_{V_1}(v) - d_{E_1}(v) - d_{E_0}(v) -2d_{V_{0,i}}(v) -1 \nonumber\\
&\geq& 2d +d_{V_1}(v) - d_{E_1}(v) - 0.2i\left(d_{V_1}(v) - d_{E_1}(v)\right) - 10^{-3}d - 2\left(0.2d_{V_0}(v)+10^{-3}d\right) -1 \nonumber\\
&=& 2d - 0.4 d_{V_0}(v) + \left(1-0.2i\right) d_{V_1}(v) - \left(1-0.2i\right)d_{E_1}(v) - 0.003d - 1 \nonumber\\
&\geq& 2d - 0.4 d_{V_0}(v) + \left(1-0.2i\right)d_{V_1}(v) - \left(1-0.2i\right)\left(0.08d_{V_1}(v)+ 10^{-3}d\right) - 0.003d - 1 \nonumber\\
&=& 2d - 0.4 d_{V_0}(v) + \left(1-0.2i\right)\cdot 0.92d_{V_1}(v) - 0.004d + 2i\cdot 10^{-4}d - 1 \nonumber\\
&=& 2d - 0.4 d_{V_0}(v) + \left(1-0.2i\right)\cdot 0.92\left(d-d_{V_0}(v)\right) - 0.004d + 2i\cdot 10^{-4}d - 1 \nonumber\\
&=& \left(2.916 - 0.184i\right)d - \left(1.32 - 0.184i\right)d_{V_0}(v)+ 2i\cdot 10^{-4}d - 1 \nonumber\\
&\geq& \left(2.916 - 0.184i\right)d - \left(1.32 - 0.184i\right) \cdot 0.0503 d + 2i\cdot 10^{-4}d - 1 \nonumber\\
&=& \left(2.849604 - 0,1745448i\right)d - 1, \label{LowerSvBoundinV0i}
\end{eqnarray}
while by~(\ref{dE0inV0}), (\ref{dE1inV0}) and (\ref{dV1}), also for every $i=0,1,2,3,4$ and each $v\in V_{0,i}$:
\begin{eqnarray}
\sigma(v) &\leq& 2d_{E_0}(v)+3\left(d-d_{E_0}(v)\right)  
= 
3d-d_{E_0}(v)  
\leq 
3d - 0.2i\left(d_{V_1}(v)-d_{E_1}(v)\right) + 10^{-3}d \nonumber\\
&\leq& 3d - 0.2i\left(d_{V_1}(v)-0.08d_{V_1}(v) - 10^{-3}d\right) + 10^{-3}d  
= 
\left(3.001 +2\cdot 10^{-4}i\right)d - 0.184i \cdot d_{V_1}(v) \nonumber\\
&\leq& \left(3.001 +2\cdot 10^{-4}i\right)d - 0.184i \left(0.95d - 0.0003d \right)  
= 
\left(3.001 - 0.1745448 i\right)d. \label{UpperSvBoundinV0i}
\end{eqnarray}
Thus by~(\ref{UpperSvBoundinV0i}) and~(\ref{LowerSvBoundinV0i}), for every $i=0,1,2,3$ and $u\in V_{0,i}$, $v\in V_{0,i+1}$:
\begin{equation}\label{V0i+1SmallerV0i}
\sigma(v) \leq \left[3.001 - 0.1745448 \left(i+1\right)\right]d = \left(2.8264552 - 0.1745448i\right)d < \sigma(u).
\end{equation}
Finally, to justify that there are no sum conflicts between vertices from $V_0$ and $V_1$, by~(\ref{V1iSmallerV1i+1}) and~(\ref{V0i+1SmallerV0i}) it is sufficient to show that sums in $V_{1,10^4}$ are smaller than sums in $V_{0,4}$.
To see that this is actually true, note that by~(\ref{SecondSv}), (\ref{DefinitionOfc1i}), (\ref{Delta2})   and~(\ref{LowerSvBoundinV0i}), for any $u\in V_{1,10^4}$ and $v\in V_{0,4}$:
\begin{eqnarray}
\sigma(u) &\leq& \left(1.0776+1.06875\cdot \left(10^4-1\right)\cdot 10^{-4}\right)d + 0.0018956d + 10^{-4}\cdot 1.0603 d \nonumber\\
&=& 2.148244755 d
<  2.1514248 d - 1
=  \left(2.849604 - 0,1745448\cdot 4\right)d - 1 
\leq \sigma(v). \nonumber
\end{eqnarray}
This finishes the proof of Theorem~\ref{123largeregTh} as the obtained weighting $\omega_2$ is thus indeed a $3$-weighting of the edges of $G$ such that there are no sum conflicts between neighbours in $G$.
\qed

\section{Concluding Remarks}

The constant $10^8$ above could still be improved, but at the cost of clarity of presentation of the proof of Theorem~\ref{123largeregTh}. Nevertheless, we were far from being able to push it down to $10^7$. 
Actually, introducing the special subgraph $G'_1$ of $G_1$, based on Corollary~\ref{QuarterDecompositionLemma}, served merely optimization purposes. These required also using only $1$'s and $3$'s as weights in $G_1$. 
However, forgetting of this direction towards optimization, focused on the lower bound for $d$, might be otherwise beneficial. Namely, using mostly $2$'s and $1$'s in $G_1$ we might assure via a similar argument as in Section~\ref{SectionProofLargeD} an arbitrarily small fraction of all edges weighted $3$ in a vertex-colouring $3$-weighting 
of any $d$-regular graph with $d$ large enough. 
Apart from this, our approach can also be relatively easily extended to graphs which are not regular, but whose minimum degree $\delta$ is slightly larger (by an arbitrary $\epsilon>0$) than half of the maximum degree $\Delta$ (and $\Delta$ is large enough)  -- this greatly improves the mentioned result from~\cite{123dense-Zhong} that the 1--2--3 Conjecture holds if $\delta>0.99985n$, where $n$ is sufficiently large order of a graph. 
We omit details here, as we believe that in fact this can still be improved towards a stronger result for general graphs for which 
it is sufficient that $\delta$ is at least a very small function of $\Delta$, of order much less than $\Delta$.
This shall however require a few extra ideas, as our approach does not directly transfer 
at this point to such a case.

\end{document}